\documentclass{amsart}
\usepackage{amssymb, amsmath}

\newcommand{\lo}{\longrightarrow}

\newcommand{\sap}{\operatorname{SAP}}
\newcommand{\wap}{\operatorname{WAP}}

\begin{document}

\title[amenability]{amenability of the algebras $R(S)$, $F(S)$ of a topological semigroup $S$}
\author[M. Amini, A.R. Medghalchi]{Massoud Amini, Alireza Medghalchi}
\address{Department of Mathematics\\Shahid Beheshti University\\ Evin, Tehran 19839, Iran\\m-amini@cc.sbu.ac.ir
\linebreak
\indent Department of Mathematics and Statistics\\ University of Saskatchewan
\\106 Wiggins Road, Saskatoon\\ Saskatchewan, Canada S7N 5E6\\mamini@math.usask.ca
\linebreak
\linebreak
\indent Department of Mathematics\\ Teacher Training University, Tehran, Iran}
\thanks{This research was partly supported by IPM}

\begin{abstract}
For a locally compact Hausdorff semigroup $S$, the $L^\infty$ representation
algebra $R(S)$ was extensively studied by Dunkl and Ramirez. The Fourier-Stieltjes
 algebra $F(S)$ of a topological semigroup was studied by Lau. The aim 
of this paper is to investigate these two algebras and study the amenability 
of them with respect to the structure of $S$. 
\end{abstract}
\maketitle

\section{Introduction.} 

Dunkl and Ramirez defined the subalgebra $R(S)$
of the algebra of weakly almost periodic functions on $S$, $\wap(S)$ [2]. 
This algebra is called the $L^\infty$-representation algebra of $S$. 
In fact $R(S)$ is the set of all functions $f(x)=\int_X(Tx)gd\mu$ $(x\in S)$ 
where $(T,X,\mu)$ is an $L^\infty$-representation of $S$ and $g\in L^1(X,\mu)$.

In [3] Lau studied the subalgebra $F(S)$ of $\wap(S)$ of a topological 
semigroup $S$ with involution. If $S$ is commutative, then $F(S)\subseteq 
R(S)$ and in particular, if $G$ is an abelian topological group, then 
$F(G)=R(G)=\widehat{M(\hat{G})}$ where $\hat{G}$ is the dual group of $G$. If 
$S$ is a topological* -
 semigroup with an identity, then $F(S)$ is the linear 
span of positive definite functions on $S$. By ([3], theorem 3.2) $F(S)$ is a 
subalgebra of $\wap(S)$. In this paper, we investigate the structure of $R(S), 
F(S)$ and $\overline{R(S)}, \overline{F(S)}$,
 the normed closures of $R(S)$ and 
$F(S)$ and show that these are left introverted subalgebras of 
$\wap(S)$ and study left (resp. right) amenability of $R(S)$, $F(S)$, 
$\overline{R(S)}$, $\overline{F(S)}$. Also we show 
that $\sap(S)\subseteq\overline{F(S)}$, where 
$\sap(S)$ is the Banach algebra of strongly almost periodic functions on $S$, 
and then show that $\overline{F(S)}$ is left amenable if and only if 
$\overline{F(S)}=\sap(S)\oplus C$, where $C$ is a closed linear subspace of 
$\overline{F(S)}$.

\section{Preliminaries} Let $S$ be a topological semigroup. Let 
$\lambda$ be a probability measure on a measureable space $X$. We know that 
$L^\infty(X,\lambda)$ is a commutative $W^*$-algebra. An 
$L^\infty$-representation $(T,X,\mu)$ of $S$ is a weak$^*$-continuous 
homomorphism $T$ of $S$ into the unit ball of $L^\infty(X,\mu)$.
\paragraph{Definition 2.1.} The representation algebra $R(S)$ is the set of 
all functions $f$ such that $f(x)=\int_XT_x(g)d\mu$, where $(T,X,\mu)$ is an 
$L^\infty$-representation of $S$ and $g\in L^1(X,\mu)$. 
\paragraph{Proposition 2.2.} $R(S)$ is a normed subalgebra of $\wap(S)$ with 
pointwise multiplication. It is closed
 under conjugation, translation invariant and 
contains constant functions. It is complete with its norm.
\paragraph{Proof.} See Theorem 1.6 and Proposition 1.4 of [2].

\vspace*{0.3cm}
Following ([1] pages 72, 77 and section 3) we have the following definitions:
\paragraph{Definition 2.3.} Let $F$ be a translation invariant subalgebra of 
$B(S)$ (the space of all bounded
 functions on $S$). $F$ is said $m$-left introverted if $T_\mu 
F\subseteq F$ $(\mu\in MM(F))$ where $T_\mu f(x)=\mu(l_xf)$, $l_x f(y)=f(xy)$ 
and $MM(F)$ the set of all multiplicative means on $F$.
\paragraph{Definition 2.4.} An admissible subalgebra of $B(S)$ is a normed 
closed, conjugate closed, translation invariant, left $m$-intrverted 
subalgebra of $B(S)$ containing constant functions. 

Following Lau we have the following definition [3].
\paragraph{Definition 2.5.} Let $S$ be a topological semigroup with 
continuous involution 
$*$. Let $M$ be a $W^*$-algebra and $M_1=\{x\in M| \|x\|\leq 1\}$. By a 
*-representation of $S$ we mean $(\omega,M)$,
where $\omega$ is a $*$-homomorphism of 
$S$ into $M_1$. 
Let $\Omega(S)$ be the set of all $\sigma$-continuous 
$*$-representations of $S$ such that $\overline{\omega(S)}^\sigma=M$. Let 
$F(S)$ be the set of all functions $f$ on $S$ such that $f=\hat{\psi}$ for 
some $\psi\in M_*$ (predual of $M$), where $\hat{\psi}=\psi\circ\omega$. Let 
$f\in F(S)$, in [3] Lau defines $\|f\|_\Omega=\inf\{\|\psi\| |\psi\in M_*, 
\hat{\psi}=f\; \text{for some} \; (\omega,M)\in\Omega(S)\}$. Note that 
$\sigma=\sigma(M,M_*)$.

\vspace*{0.3cm}
Following [3] we have the following proposition. 
\paragraph{Proposition 2.6.} $F(S)$ is a commutative 
subalgebra of $\wap(S)$ which is conjugate closed, translation invariant, 
containing the constant functions. Also, if $f\in F(S)$, then $f^*\in
F(S)$, where $f^*(s)=\overline{f(s^*)}$. Furthermore, $\|\cdot\|_\Omega$ is a 
norm on $F(S)$ and $(F(S),\|\cdot\|_\Omega)$ is a commutative normed algebra 
with unit.
\paragraph{Proof.} This is theorem 3.2 of [3] where its proof was refered to 
[2].
\section{Amenability of $F(S), R(S), \overline{F(S)}, \overline{R(S)}$.}
In this section we assume that $S$ is a topological semigroup with a 
continuous involution.

In the following theorem $K(X), K(Y)$ are the minimal ideals of $X,Y$, 
respectively.
\paragraph{Theorem 3.1.} Let $Y=\sigma(F(S))$ (the spectrum of $F(S)$) 
which is a semigroup by proposition 5.4 in [3]. Then $F(S)$ 
is amenable if and only if $K(Y)$ is a topological group. 
\paragraph{Proof.} Assume $F(S)$ is amenable and let 
 $\overline{F(S)}$ be the norm closure of $F(S)$ in 
$\wap(S)$. Then clearly $\overline{F(S)}\subseteq\wap(S)$. By (theorem 3.2 
[3]), $\|f\|_\infty\leq\|f\|_\Omega$ $(f\in F(S))$. So the embedding function 
$i:F(S)\lo \overline{F(S)}$ is continuous. Therefore by (theorem 1.9. [5]), 
$\overline{F(S)}$ is amenable. Now again by (theorem 3.2, [1]) $F(S)$ is a 
translation invariant, conjugate closed subalgebra of $\wap(S)$ and contains 
constant functions, hence $\overline{F(S)}$ is norm closed, trnaslation 
invariant, conjugate closed subalgebra of $\wap(S)$ containing constant 
functions. Therefore by (corollary 4.2.7, [1]) $\overline{F(S)}$ is 
introverted. So, it is an $m$-admissible subalgebra of $\wap(S)$. Hence, by 
(theorem 4.2.14, [1]) $K(X)$ is a topological 
group, where $X$ is $S^{\overline{F(S)}}$, the maximal ideal 
space of $\overline{F(S)}$. Now let $\pi=i^*:\overline{F(S)}^*\lo 
F(S)^*$ is the adjoint mapping of $i$. Clearly $\pi$ is norm decreasing and 
onto and $\pi(X)=Y$. Thus $\pi(K(X))=K(Y)$ is a topological group. Conversely, 
if $K(Y)$ is a topological group, then $K(X)$ is a topological group, so by 
(4.2.14, [1]) $\overline{F(S)}$ is amenable, and so is $F(S)$.

\vspace*{0.3cm}
We have another consequence of (theorem 4.2.14, [1]). In fact, 
$\overline{F(S)}$ if left (right) amenable if and only if $K(X)$ is a minimal 
right (left) ideal of $X$. 

Now, we extend the above theorem for left (resp. right) amenability of $F(S)$.
\paragraph{Theorem 3.2.} $F(S)$ is left (right) amenable if and only if $K(Y)$ 
is minimal right (left) ideal of $Y$. 
\paragraph{Proof.} Let $F(S)$ be left
 amenable and $f\in\overline{F(S)}$. Then 
there is a sequence $\{g_n\}$ in $F(S)$ such that 
$g_n\overset{\|\cdot\|_\infty}{\lo} f$. Clearly $\{m(g_n)\}$ is a Cauchy 
sequence in $\overline{F(S)}$, so it is convergent. We define 
$\bar{m}(f)=\lim_{n\rightarrow\infty} m(g_n)$. Since $\bar{m}$ is positive, 
$\bar{m}(1)=1$ and 
$\bar{m}(l_xf)=\lim_{n\rightarrow\infty}m(_xg_n)=\lim_{n\rightarrow \infty} 
m(g_n)=\bar{m}(f)$, it follows that $\bar{m}$ is a left invariant mean on 
$\overline{F(S)}$. Hence, $\overline{F(S)}$ is left amenable. By (theorem 
4.2.14, [1]) $K(X)$ is a minimal right ideal of $X$. Now, we have 
$\pi(K(X))=K(Y)$, therefore by (corollary 1.3.16, [1]) $K(Y)$ is a minimal 
right  ideal of $Y$. Conversely, if $K(Y)$ is a minimal right ideal of $Y$, 
then by the same argument
 $K(X)$ is a minimal right ideal of $X$ and therefore, 
again by (theorem 4.2.14, [1]) $\overline{F(S)}$ is left amenable. Thus $F(S)$ 
is left amenable.
The right version is proved in a similar way. 

Next we show that $F(S)$ is an $F$-algebra in the sense of Lau [4].
\paragraph{Proposition 3.3.} If $S$ is a unital topological semigroup with 
continuous involution, then $F(S)$ is an $F$-algebra. 
\paragraph{Proof.} Each element of $F(S)$ is of the form 
$\hat{\psi}=\psi\circ\omega_\Omega$, for some $\psi\in(M_\Omega)_*$, where 
$(\omega_\Omega,M_\Omega)$, is the universal representation of $S$ [3]. Take 
$\psi_1,\psi_2\in(M_\Omega)_*$, then by the Gelfand-Naimark-Segal 
construction, there are (unit) vectors $\xi_i,\eta_i\in H_\Omega$ such that 
$\psi_i(x)=<x\xi_i,\eta_i>$ $(x\in M_\Omega, i=1,2)$. Then 
$$\hat{\psi}_1\hat{\psi}_2(s)=\hat{\psi}_1(S)\hat{\psi}_2(s)=<\omega_\Omega(s) 
\xi_1,\eta_1><\omega_\Omega(s)\xi_2,\eta_2>= <\omega_\Omega(s) 
\xi_1\otimes\xi_2, \eta_1\otimes\eta_2> \; (s\in S).$$
Hence, if $1\in W^*_\Omega(S)$ is the identity element (see [3] for the 
notation), then 
$$<1,\hat{\psi}_1\hat{\psi}_2>=<\xi_1\otimes\xi_2, \eta_1\otimes\eta_2>= 
<\xi_1,\eta_1><\xi_2,\eta_2>=<1,\hat{\psi}_1><1,\hat{\psi}_2>,$$
and we are done.
\paragraph{Corollary 3.4.} If $S$ is a unital topological semigroup with 
continuous involution, then $W^*_\Omega(S)$ has a topological left invariant 
mean.
\paragraph{Proof.} This follows from above proposition and (theorem 4.1, [4]). 
Note that $F(S)$ is always left amenable in the sense of [4]. 
\paragraph{Remark 3.5.} (a) If $S$ is an idempotent commutative topological 
semigroup with involution $s=s^*$, then by (3.3(c), [3]) $F(S)=R(S)$. 
Therefore in this special case the above results hold for $R(S)$ too.

(b) By Proposition 2.2, $R(S)$ is right translation invariant, conjugate 
closed, containing constant functions. Therefore when $S$ is commutative,
$\overline{R(S)}\subseteq\wap(S)$ is norm closed, translation invariant, 
conjugate closed  subalgebra of $\wap(S)$ containing constant functions.
therefore $\overline{R(S)}$ is an $m$-admissible subalgebra of $\wap(S)$ 
(theorem 4.2.14, [1]). Hence, the results of theorems 3.1 and 3.2 hold for 
$\overline{R(S)}$.

(c) If $G$ is a non-compact locally compact abelian group, then $\wap(G)\neq 
\overline{R(G)}=\overline{\widehat{M(\hat{G})}}$ (see 52.10, [2]), or 
even if $S$ is a non compact locally compact subsemigroup of a locally compact 
group, then $\wap(S)\neq\overline{R(S)}$ (see (5.2.12, [2]).

(d) Example 4.2 of [3] also shows that $\overline{F(S)}
\neq\overline{R(S)}$ may happen.

(e) When $S$ is an abelian semigroup with involution, we have another 
interesting result. By (4.3.8, [1]) $\sap(S)$ (the Banach space of all strongly
 almost periodic 
functions on $S$) is the closed linear span of characters. Therefore, by 
(3.3.(e), [1]) $\sap(S)\subseteq F(S)\subseteq\wap(S)$.

Now by (theorem 4.3.12, [1]) $\overline{F(S)}$ is left amenable if and only if 
$\overline{F(S)}=\sap(S)\oplus C$ where $C$ is a left translation 
invariant, closed linear subspace of $\overline{F(S)}$, if and only if 
there exists a projection of $\overline{F(S)}$ onto $\sap(S)$ that 
commutes with left translations. 
Since, the amenability of $F(S)$ and $\overline{F(S)}$ are 
equivalent, the amenability of $F(S)$ is equivalent to each of the above 
conditions. 

(f) If $S$ is a topological group, then $Y=\sigma(B(S))$ is a $\sigma$-weakly 
compact, self adjoint subset of the unit ball of $W^*(S)$ and $Y\cap 
W^*(S)_u\simeq S$ (see [6] for details and notations). If $X$ is 
$\overline{B(S)}$-compactification of $S$ and $G_1=K(Y)$, $G_2=K(X)$, then by 
theorem 3.1. and (3.4, [3]), the Banach algebras $A(S)$ and $B(S)/A(S)$ are 
both amenable iff any of $G_1$ or $G_2$ is a topological group. This is 
interesting in the light of the fact that we don't know exactly when $A(S)$ is 
amenable.

Next we show that remark (e) indeed applies to any topological $*$-semigroup.
\paragraph{Lemma 3.6.} Let $S$ be a topological $*$-semigroup, then 
$$\sap(S)\subseteq\overline{F(S)}\subseteq\wap(S).$$
\paragraph{Proof.} Let $\{\pi,H\}$ be a (finite dimensional) unitary 
representation of $S$, $\xi,\eta\in H$, and $u(s)=<\pi(s)\xi,\eta>$ $(s\in 
S)$. Put $M=<\pi(S)>^{-\sigma}\subseteq B(H)$ and 
$\alpha=(\pi,M)\in\Omega(S)$. Consider $\psi(x)=<x\xi,\eta>$ $(x\in M)$. Then 
obviousely $\psi\in M_*$. Indeed $\psi$ is $w^*$-continuous on $B(H)$: if 
$\{T_\alpha\}\subseteq B(H)$ and $T_\alpha\overset{w^*}{\lo} T$ in $B(H)$, 
then considering the rank one operator $\xi\otimes\eta$ in $B(H)_*$, we have
$$<T_\alpha\xi,\eta>=<T_\alpha,\xi\otimes\eta>\underset{\alpha}{\lo} 
<T,\xi\otimes\eta>=<T\xi,\eta>$$
i.e. $\psi(T_\alpha)\lo\psi(T)$. Now we have $\hat{\psi}=\psi\circ\pi=u$, 
hence $u\in F(S)$. Therefore $\sap(S)\subseteq\overline{F(S)}$. The other 
inclusion is contained in [3].

Now the proof of the following theorem is exactly as in remark (e). 
\paragraph{Theorem 3.7.} Let $S$ be a topological  semigroup with 
continuous involution. Then 
the following are equivalent:
\begin{itemize}
\item[a)] $R(S)$ is amenable. 
\item[b)] $F(S)$ is amenable. 
\item[c)] $\overline{R(S)}$ is amenable. 
\item[d)] $\overline{F(S)}$ is amenable. 
\item[e)] $\overline{F(S)}=\sap(S)\oplus C$, where $C$ is a translation 
invariant,  closed linear subspace of $\overline{F(S)}$. 
\item[f)] There exists a projection of $F(S)$ onto $\sap(S)$ that 
commutes with translations. 
\end{itemize}
%%%%%%%%%%%


\newpage
\begin{thebibliography}{99}
\bibitem{} Berglund, John F., Junghenn, Hugo, D., Milnes, P.,  Analysis on 
semigroups, John Wiley \& Sons, 1989.
\bibitem{} Dunkl, Charles, F., Ramirez, Donald, E.,
Representation of commutative semitopological semigroup, 
 Lecture Notes in Mathematics, 435(1975) Springer-Verlag.
\bibitem{} Lau, A. T., The Fourier Stieltjes algebra of a topological 
semigroup with involution, Pacific J.  Math. vol. 77, No. 1 (1978), 165-181.
\bibitem{} Lau, A. T., Analysis on a class of Banach algebras with 
applications to harmonic analysis on locally compact groups and semigroups, 
Fundamenta Math. vol. 68 (1983), 161-175.
\bibitem{} Pier, Jean. Paul, Amenable Banach algebras, John Wiley, 1988.
\bibitem{} Walter, Martin E., $W^*$-algebras and nonabelian Harmonic analysis, 
J. Func. Anal. vol. 11 (1972), 17-38.
\end{thebibliography}
\end{document}